\documentclass[11pt]{amsart}
\usepackage{amsmath,amssymb,amsfonts,url}
\numberwithin{equation}{section}

\newtheorem{thm}{Theorem}[section]
\newtheorem{lem}{Lemma}[section]
\newtheorem{rem}{Remark}[section]
\newtheorem{prop}{Proposition}[section]
\newtheorem{cor}{Corollary}[section]

\begin{document}
\title[Monge-Amp\`{e}re equation]{Monge-Amp\`{e}re equation on exterior domains} 
\author{Jiguang Bao}
\address{School of Mathematical Sciences\\
Beijing Normal University, Laboratory of Mathematics and Complex Systems\\
Ministry of Education \\
Beijing 100875, China}
\email{jgbao@bnu.edu.cn}

\author{Haigang Li}
\address{School of Mathematical Sciences\\
Beijing Normal University, Laboratory of Mathematics and Complex Systems\\
Ministry of Education \\
Beijing 100875, China}
\email{hgli@bnu.edu.cn}

\author{Lei Zhang}
\address{Department of Mathematics\\
        University of Florida\\
        358 Little Hall P.O.Box 118105\\
        Gainesville FL 32611-8105}
\email{leizhang@ufl.edu}
\thanks{Bao is supported by NSFC (11071020) and SRFDPHE (20100003110003), Li is supported by NSFC (11071020)(11126038)(11201029) and SRFDPHE (20100003120005), Zhang is supported in part by NSF Grant 1027628}

\date{\today}

%%%%%%%%%%%%%%%%%%%%%%%%%%%%%%%%%%%%%%%%%%%%%
\begin{abstract} We consider the Monge-Amp\`ere equation $\det(D^2u)=f$ where $f$ is a positive function in $\mathbb R^n$ and $f=1+O(|x|^{-\beta})$ for some $\beta>2$ at infinity. If the equation is globally defined on $\mathbb R^n$ we classify the asymptotic behavior of solutions at infinity. If the equation is defined outside a convex bounded set we solve the corresponding exterior Dirichlet problem. Finally we prove for $n\ge 3$ the existence of global solutions with prescribed asymptotic behavior at infinity. The assumption $\beta>2$ is sharp for all the results in this article.

\end{abstract}
%%%%%%%%%%%%%%%%%%%%%%%%%%%%%%%%%%%%%%%%%%%%%

\maketitle

\section{Introduction}

It is well known that Monge-Amp\`ere equations are a class of important fully nonlinear equations profoundly related to many fields of analysis and geometry. In the past few decades many significant contributions have been made on various aspects of Monge-Amp\`ere equations. In particular, the Dirichlet problem
$$\left\{\begin{array}{ll}
det(D^2u)=f, \quad \mbox{in}\quad D, \\
u=\phi \quad \mbox{on }\quad \partial D
\end{array}
\right.
$$
on a convex, bounded domain $D$ is completely understood through the works of Aleksandrov \cite{alek1}, Bakelman \cite{b1}, Nirenberg \cite{n2}, Calabi \cite{calabi}, Pogorelov \cite{p1,p5,p8}, Cheng-Yau \cite{cheng-yau1}, Caffarelli-Nirenberg-Spruck \cite{cns}, Caffarelli \cite{caf-ann2}, Krylov \cite{k2}, Jian-Wang \cite{jian-wang}, Huang \cite{huang}, Trudinger-Wang \cite{tw7}, Urbas \cite{urbas},Savin \cite{savin1,savin2}, Philippis-Figalli \cite{figalli}  and the references therein.
Corresponding to the traditional Dirichlet problem mentioned above, there is an exterior Dirichlet problem which seeks to solve the Monge-Amp\`ere equation outside a convex set. More specifically, let $D$ be a smooth, bounded and strictly convex subset of $\mathbb R^n$ and let $\phi\in C^2(\partial D)$, the exterior Dirichlet problem is to find $u$ to verify
\begin{equation}\label{extdir}
\left\{\begin{array}{ll}
det(D^2u)=f(x),\quad x\in \mathbb R^n\setminus \overline{D}, \\
u\in C^0(\mathbb R^n\setminus D) \mbox{ is a locally convex viscosity solution},\\
u=\phi,\quad \mbox{on }\quad \partial D.
\end{array}
\right.
\end{equation}

If $f\equiv 1$ and $n\ge 3$, Caffarelli and Li \cite{caf-li1} proved that any solution $u$ of (\ref{extdir}) is very close to a parabola near infinity. They solved the exterior Dirichlet problem assuming that $u$ equals $\phi$ on $\partial D$ and has a prescribed asymptotic behavior at infinity. For $f\equiv 1$ and $n=2$,
Ferrer-Mart\'inez-Mil\'an \cite{FMM1,FMM2} used a method of complex analysis to prove that any solution $u$ of (\ref{extdir}) is very close to a parabola plus a logarithmic function at infinity (see also Delano\"e \cite{delanoe}). Recently the first two authors \cite{bao} solved the
exterior Dirichlet problem for $f\equiv 1$ and $n=2$. In the first part of this article we solve the exterior Dirichlet problem assuming that $f$ is a perturbation of $1$ near infinity:
\begin{eqnarray*}
(FA)&:& f\in C^0(\mathbb R^n),\,\, 0<\inf_{\mathbb R^n} f\le \sup_{\mathbb R^n} f<\infty. \\
&& D^3f\, \mbox{exists outside a compact subset of } \mathbb R^n,\\
&& \exists\beta>2 \mbox{ such that }
 \lim_{|x|\to \infty} |x|^{\beta+|\alpha |} |D^{\alpha}(f(x)-1)|<\infty, \,\, |\alpha |=0,1,2,3.
\end{eqnarray*}

Let $\mathbb M^{n\times n}$ be the set of the real valued, $n\times n$ matrices and
$$\mathcal{A}:=\{A\in \mathbb M^{n\times n}:\quad A \mbox{ is symmetric, positive definite and }\,\,\, det(A)=1 \, \}. $$
Our first main theorem is

\begin{thm}\label{thm5} Let $D$ be a strictly convex, smooth and bounded set, $\phi\in C^2(\partial D)$ and $f$ satisfy (FA).
If $n\ge 3$, then for any $b\in \mathbb R^n$, $A\in \mathcal{A}$, there exists $c_*(n, D, \phi, b, A, f)$ such that for any $c>c_*$, there exists a unique $u$ to (\ref{extdir}) that satisfies
\begin{equation}\label{11may5e1}
\limsup_{|x|\to \infty} |x|^{\min\{\beta,n\}-2+k}|D^k(u(x)-(\frac 12x'Ax+b\cdot x+c))|<\infty
\end{equation}
for $k=0,1,2,3,4$.
If $n=2$, then for any $b\in \mathbb R^2$, $A\in \mathcal{A}$, there exists $d^*\in \mathbb R$ depending only on $A,b,\phi,f,D$ such that for all $d>d^*$, there exists a unique $u$ to (\ref{extdir}) that satisfies
\begin{equation}\label{12jan18e1}
\limsup_{|x|\to \infty}|x|^{k+\sigma}\bigg |D^k(u(x)-(\frac 12 x'Ax+b\cdot x+d\log \sqrt{x'Ax}+c_d)) \bigg |<\infty
\end{equation}
for $k=0,1,2,3,4$ and $\sigma\in (0,\min\{\beta-2,2\})$. $c_d\in \mathbb R$ is uniquely determined by $D,\phi,d,f,A,b$.
\end{thm}

\medskip

The Dirichlet problem on exterior domains is closely related to asymptotic behavior of solutions defined on entire $\mathbb R^n$. The classical theorem of J\"orgens \cite{jorgens}, Calabi \cite{calabi} and Pogorelov \cite{pogorelov1} states that any convex classical solution of $det(D^2u)=1$ on $\mathbb R^n$ must be a quadratic polynomial. See
Cheng-Yau\cite{cheng-yau2}, Caffarelli \cite{caf-ann2} and Jost-Xin \cite{jost-xin} for different proofs and extensions. Caffarelli-Li\cite{caf-li1} extended this result by considering
\begin{equation}\label{12feb27e1}
det(D^2 u)=f \quad \mathbb R^n
\end{equation}
where $f$ is a positive continuous function and is not equal to $1$ only on a bounded set. They proved that for $n\ge 3$, the convex viscosity solution $u$ is very close to quadratic polynomial at infinity and for $n=2$, $u$ is very close to a quadratic polynomial plus a logarithmic term asymptotically. In a subsequent work \cite{caf-li2} Caffarelli-Li proved that if $f$ is periodic, then $u$ must be a perturbation of a quadratic function.

The second main result of the paper is to extend the Caffarelli-Li results on global solutions in \cite{caf-li1}:

\begin{thm}\label{thm1} Let $u\in C^0(\mathbb R^n)$ be a convex viscosity solution to (\ref{12feb27e1})
where $f$ satisfies $(FA)$.
If $n\ge 3$,
then there exist $c\in \mathbb R$, $b\in \mathbb R^n$ and $A\in \mathcal{A}$ such that (\ref{11may5e1}) holds. If $n=2$ then
there exist $c\in \mathbb R$, $b\in \mathbb R^2$, $A\in \mathcal{A}$ such that (\ref{12jan18e1}) holds for
$d=\displaystyle{\frac 1{2\pi}\int_{\mathbb R^2}(f-1)}$ and $\sigma\in (0,\min\{\beta-2,2\})$.
\end{thm}

\begin{cor}\label{cor1}
Let $D$ be a bounded, open and convex subset of $\mathbb R^n$ and let $u\in C^0(\mathbb R^n\setminus \bar D)$ be a locally convex viscosity solution to
\begin{equation}\label{12apr5e1}
det(D^2u)=f, \quad \mbox{in }\quad \mathbb R^n\setminus \bar D
\end{equation}
where $f$ satisfies (FA). Then for $n\ge 3$, there exist $c\in \mathbb R$, $b\in \mathbb R^n$ and $A\in \mathcal{A}$ such that
(\ref{11may5e1}) holds. For $n=2$, there exist $A\in \mathcal{A}$, $b\in \mathbb R^n$ and $c, d\in \mathbb R$ such that for $k=0,..,4$
$$\limsup_{|x|\to \infty}|x|^{k+\sigma}\bigg |D^k(u(x)-(\frac 12 x'Ax+b\cdot x+d\log \sqrt{x'Ax}+c)) \bigg |<\infty $$
holds for $\sigma\in (0,\min\{\beta-2,2\})$.
\end{cor}

As is well known the Monge-Amp\`ere equation $det(D^2u)=f$ is closely related to the Minkowski problems, the Plateau type problems, mass transfer problems, and affine geometry, etc. In many of these applications $f$ is not a constant. The readers may see the survey paper of Trudinger-Wang \cite{tw1} for more description and applications. The importance of $f$ not identical to $1$ is also mentioned by Calabi in \cite{calabi}.

Next we consider the globally defined equation (\ref{12feb27e1})
and the existence of global solutions with prescribed asymptotic behavior at infinity.

\begin{thm}\label{thm3} Suppose $f$ satisfies (FA). Then for any $A\in \mathcal{A}$, $b\in \mathbb R^n$ and $c\in \mathbb R$, if $n\ge 3$
there exists a unique convex viscosity solution $u$ to (\ref{12feb27e1}) such that (\ref{11may5e1}) holds.
\end{thm}

\medskip

The following example shows that the decay rate assumption $\beta>2$ in (FA) is sharp in all the theorems. Let $f$ be a radial, smooth, positive function such that $f(r)\equiv 1$ for $r\in [0,1]$ and $f(r)=1+r^{-2}$ for $r>2$. Let
$$u(r)=n^{\frac 1n}\int_0^r\bigg (\int_0^s t^{n-1}f(t)dt\bigg )^{\frac 1n}ds,\quad r=|x|. $$
It is easy to check that $det(D^2u)=f$ in $\mathbb R^n$. Moreover for $n\ge 3$, $u(x)=\frac 12|x|^2+O(\log |x|)$ at infinity. For $n=2$,
$u(x)=\frac 12|x|^2+O((\log |x|)^2)$ at infinity.

Corresponding to the results in this paper we make the following two conjectures. First
we think the analogue of Theorem \ref{thm3} for $n=2$ should also hold.

\emph{Conjecture 1: Let $n=2$ and $f$ satisfy (FA),
then there exists a unique convex viscosity solution $u$ to (\ref{12feb27e1}) such that
(\ref{12jan18e1}) holds for $d=\frac{1}{2\pi}\int_{\mathbb R^2}(f-1)$ and $\sigma\in (0,\min\{\beta-2,2\})$.}

\medskip

\emph{Conjecture 2: The $d^*$ in Theorem \ref{thm5} is $\frac 1{2\pi}\int_{\mathbb R^2\setminus D}(f-1)-\frac{1}{2\pi}\mbox{area}(D)$. }

\medskip

These two conjectures are closely related in a way that if conjecture one is proved, then conjecture two follows by the same argument in the proof of Theorem \ref{thm5}.

The organization of this paper is as follows: First we establish a useful proposition in section two, which will be used in the proof of all theorems. Then in section three we prove Theorem \ref{thm5} using Perron's method. Theorem \ref{thm3} and Theorem \ref{thm1} are proved in section four and section five, respectively. In the appendix we cite the interior estimates of Caffarelli and Jian-Wang. The proof of all the theorems in this article relies on previous works of Caffarelli \cite{caf-ann1,caf-ann2},
Jian-Wang \cite{jian-wang} and Caffarelli-Li \cite{caf-li1}. For example, Caffarelli-Li \cite{caf-li1} made it clear that for exterior Dirichlet problems, convex viscosity solutions are strictly convex. On the other hand for Monge-Amp\`ere equations on convex domains, Pogorelov has a well known example of a not-strictly-convex solution. Besides this, we also use the Alexandrov estimates, the interior estimate of Caffarelli \cite{caf-ann2} and  Jian-Wang \cite{jian-wang} in an essential way.

\section{A useful proposition}

Throughout the article we use $B_r(x)$ to denote the ball centered at $x$ with radius $r$ and $B_r$ to denote the ball of radius $r$ centered at $0$.

The following proposition will be used in the proof of all theorems.
\begin{prop}\label{caf-li-prop}
Let $R_0>0$ be a positive number, $v\in C^0(\mathbb R^n\setminus \bar B_{R_0})$ be a convex viscosity solution of
$$det(D^2v)=f_v \quad \mathbb R^n\setminus \bar B_{R_0}$$
where $f_v\in C^3(\mathbb R^n)\setminus \bar B_{R_0}$ satisfies
$$\frac 1{c_0}\le f_v(x)\le c_0, \quad x\in \mathbb R^n \setminus B_{R_0}$$
and
\begin{equation}\label{12feb23e1}
|D^k(f_v(x)-1)|\le c_0|x|^{-\beta-k},\quad |x|>R_0, \quad k=0,1,2,3.
\end{equation}
Suppose there exists $\epsilon>0$ such that
\begin{equation}\label{12feb27e2}
|v(x)-\frac 12|x|^2|\le c_1|x|^{2-\epsilon},\quad |x|\ge R_0
\end{equation}
then for $n\ge 3$, there exist $b\in \mathbb R^n$, $c\in \mathbb R$ and $C(n,R_0,\epsilon,\beta,c_0,c_1)$
such that,
\begin{eqnarray}\label{12mar7e10}
&&|D^k(v(x)-\frac 12|x|^2-b\cdot x-c)|\\
&\le &C/|x|^{\min\{\beta,n\}-2+k}, \,\, |x|>R_1,\,\, k=0,1,2,3,4;\nonumber
\end{eqnarray}
where $R_1(n,R_0,\epsilon,\beta,c_0,c_1)>R_0$ depends only on $n,R_0,\epsilon,\beta,c_0$ and $c_1$.
For $n=2$, there exist $b\in \mathbb R^2$, $d,c\in \mathbb R$
such that for all $\sigma\in (0,\min\{\beta-2,2\})$
\begin{eqnarray} \label{12mar7e11}
&&|D^k(v(x)-\frac 12|x|^2-b\cdot x-d\log |x|-c)|\\
&\le &\frac{C(\epsilon,R_0,\beta,c_0,c_1)}{|x|^{\sigma+k}},\,\, |x|>R_1, \,\, k=0,1,2,3,4 \nonumber
\end{eqnarray}
where $R_1>R_0$ depends only on $\epsilon,R_0,\beta,c_0$ and $c_1$.
\end{prop}

\begin{rem} $\epsilon$ may be greater than or equal to $2$ in Proposition \ref{caf-li-prop}.
\end{rem}

\noindent{\bf Proof of Proposition \ref{caf-li-prop}:}
Proposition \ref{caf-li-prop} is proved in \cite{caf-li1} for the case that $f\equiv 1$ outside a compact subset of $\mathbb R^n$. For this more general case, Theorem \ref{caf-jian-wang} in the appendix (A theorem of Caffarelli, Jian-Wang) and Schauder estimates play a central role.
First we establish a lemma that holds for all dimensions $n\ge 2$.

\begin{lem}\label{12feb27lem1} Under the assumption of Proposition \ref{caf-li-prop}, let
$$w(x)=v(x)-\frac 12|x|^2, $$
then there exist $C(n,R_0,\epsilon,c_0,c_1,\beta)>0$ and $R_1(n,R_0,\epsilon,c_0,c_1,\beta)>R_0$ such that
for any $\alpha\in (0,1)$
\begin{equation}\label{12feb22e4}
\left\{\begin{array}{ll}
|D^kw(y)|\le C|y|^{2-k-\epsilon_{\beta}}, \quad k=0,1,2,3,4,\quad |y|>R_1\\
\\
\frac{|D^4w(y_1)-D^4w(y_2)|}{|y_1-y_2|^{\alpha}}\le C|y_1|^{-2-\epsilon_{\beta}-\alpha},\quad |y_1|>R_1,\,\, y_2\in B_{\frac{|y_1|}{2}}(y_1)
\end{array}
\right.
\end{equation}
where $\epsilon_{\beta}=\min\{\epsilon,\beta\}$.
\end{lem}

\noindent{\bf Proof of Lemma \ref{12feb27lem1}:}

For $|x|=R>2R_0$, let
$$v_R(y)=(\frac 4R)^2v(x+\frac{R}4y),\quad |y|\le 2, $$
and
$$w_R(y)=(\frac 4R)^2w(x+\frac{R}4y),\quad |y|\le 2.  $$
By (\ref{12feb27e2}) we have
\begin{equation}\label{12feb22e10}
\|v_R\|_{L^{\infty}(B_2)}\le C, \quad \|w_R\|_{L^{\infty}(B_2)}\le CR^{-\epsilon}
\end{equation}
and
$$v_R(y)-(\frac 12|y|^2+\frac 4Rx\cdot y+\frac{8}{R^2}|x|^2)=O(R^{-\epsilon}),\quad B_2. $$
Let $\bar v_R(y)=v_R(y)-\frac 4Rx\cdot y-\frac{8}{R^2}|x|^2$, clearly $D^2\bar v_R=D^2v_R$. If $R>R_1$ with $R_1$ sufficiently large,
the set
$$\Omega_{1,v}=\{y\in B_2;\quad \bar v_R(y)\le 1 \quad \} $$
is between $B_{1.2}$ and $B_n$.
The equation for $\bar v_R$ is
\begin{equation}\label{12feb23e2}
det(D^2 \bar v_R(y))=f_{1,R}(y):=f_v(x+\frac R4y), \quad \mbox{on } B_2.
\end{equation}
Immediately from  (\ref{12feb23e1}) we have, for any $\alpha\in (0,1)$
\begin{equation}\label{12feb22e11}
\|f_{1,R}-1\|_{L^{\infty}(B_2)}+\|Df_{1,R}\|_{C^{\alpha}(B_2)}+\|D^2f_{1,R}\|_{C^{\alpha}(B_2)}\le CR^{-\beta}.
\end{equation}
Applying Theorem \ref{caf-jian-wang} on $\Omega_{1,v}$
$$\|D^2v_R\|_{C^{\alpha}(B_{1.1})}=\|D^2\bar v_R\|_{C^{\alpha}(B_{1.1})}\le C. $$
Using (\ref{12feb23e2}) and (\ref{12feb22e11}) we have
\begin{equation}\label{11may9e3}
\frac{I}C\le D^2v_R\le CI \quad \mbox{ on }B_{1.1}
\end{equation}
for some $C$ independent of $R$. Rewrite (\ref{12feb23e2}) as
$$a_{ij}^R \partial_{ij}v_R=f_{1,R},\quad B_2 $$
where $a_{ij}^R=cof_{ij}(D^2 v_R)$. Clearly by(\ref{11may9e3}) $a_{ij}^R$ is uniformly elliptic and
$$\|a_{ij}^R\|_{C^{\alpha}(\bar B_{1.1})}\le C. $$
By Schauder estimates
\begin{equation}\label{11may9e2}
 \|v_R\|_{C^{2,\alpha}(\bar B_{1})}\le C(\|v_R\|_{L^{\infty}(\bar B_{1.1})}+\|f_{1,R}\|_{C^{\alpha}(\bar B_2)})\le C.
\end{equation}

For any $e\in \mathbb S^{n-1}$, apply $\partial_e$ to both sides of (\ref{12feb23e2})
\begin{equation}\label{12feb22e13}
a^R_{ij}\partial_{ij}(\partial_e v_R)=\partial_e f_{1,R}.
\end{equation}
Since $a^R_{ij}$, $\partial_e v_R$ and $\partial_ef_{1,R}$ are bounded in $C^{\alpha}$ norm, we have
\begin{equation}\label{12feb22e15}
\|v_R\|_{C^{3,\alpha}(\bar B_{1})}\le C,
\end{equation}
which implies
\begin{equation}\label{12feb22e16}
 \|a^R_{ij}\|_{C^{1,\alpha}(\bar B_{1})}\le C.
 \end{equation}

The difference between (\ref{12feb23e2}) (with $\bar v_R$ replaced by $v_R$)
and $det(I)=1$ gives
\begin{equation}\label{11may9e8}
\tilde a_{ij}\partial_{ij}w_R=f_{1,R}(y)-1=O(R^{-\beta})
\end{equation}
where
$\tilde a_{ij}(y)=\int_0^1cof_{ij}(I+tD^2w_R(y))dt$.
By (\ref{11may9e3}) and (\ref{12feb22e15})
$$\frac{I}{C}\le \tilde a_{ij} \le CI,\,\,\mbox{on}\,\, B_{1.1},\quad \|\tilde a_{ij}\|_{C^{1,\alpha}(\bar B_{1})}\le C. $$
Thus Schauder's estimate gives
\begin{equation}\label{11may9e3a}
\|w_R\|_{C^{2,\alpha}(B_1)}\le C(\|w_R\|_{L^{\infty}(\bar B_{1.1})}+\|f_{1,R}-1\|_{C^{\alpha}(\bar B_{1})})\le CR^{-\epsilon_{\beta}}.
\end{equation}
Going back to (\ref{12feb22e13}) and rewriting it as
$$a_{ij}^R\partial_{ij}(\partial_e w_R)=\partial_e f_{1,R}. $$
We obtain, by Schauder's estimate,
\begin{equation}\label{12mar6e20}
\|w_R\|_{C^{3,\alpha}(\bar B_{1/2})}\le C(\|w_R\|_{L^{\infty}(\bar B_{3/4})}+|Df_{1,R}\|_{C^{\alpha}(\bar B_{3/4})})\le CR^{-\epsilon_{\beta}}.
\end{equation}
Since $\partial_{ije}w_R=\partial_{ije}v_R$, we also have
$$\|D^3v_R\|_{C^{\alpha}(\bar B_{1/2})}\le CR^{-\epsilon_{\beta}}.$$
By differentiating on (\ref{12feb22e13}) with respect to any $e_1\in \mathbb S^{n-1}$ we have
$$a_{ij}^R\partial_{ij}(\partial_{ee_1}w_R)=\partial_{ee_1}f_{1,R}-\partial_{e_1}a_{ij}^R\partial_{ij}\partial_ew_R. $$
(\ref{12mar6e20}) gives
$$\|\partial_{e_1}a_{ij}^R\partial_{ij}\partial_ew_R\|_{C^{\alpha}(\bar B_{1/2})}\le CR^{-2\epsilon_{\beta}}. $$
Using (\ref{12feb22e11}),$\|\partial_{ee_1}f_{1,R}\|_{C^{\alpha}(\bar B_1)}\le CR^{-\beta}$
and Schauder's estimate we have
\begin{equation}\label{12feb22e6}
\|w_R\|_{C^{4,\alpha}(\bar B_{1/4})}\le CR^{-\epsilon_{\beta}},
\end{equation}
which implies (\ref{12feb22e4}). Lemma \ref{12feb27lem1} is established. $\Box$

\medskip

Next we prove a lemma that improves the estimates in Lemma \ref{12feb27lem1}.
\begin{lem}\label{12feb27lem2} Under the same assumptions of Lemma \ref{12feb27lem1} and let $R_1$ be the large constant determined in the proof of
Lemma \ref{12feb27lem2}.  If in addition $2\epsilon<1$, then for $n\ge 3$
$$\left\{\begin{array}{ll}
|D^jw(x)|\le C|x|^{2-2\epsilon-j},\quad |x|>2R_1,\quad j=0,1,2,3,4 \\
\\
\frac{|D^4w(y_1)-D^4w(y_2)|}{|y_1-y_2|^{\alpha}}\le C|y_1|^{-2-2\epsilon-\alpha},\quad |y_1|>2R_1,\,\, y_2\in B_{|y_1|/2}(y_1)
\end{array}
\right.
$$
where $\alpha\in (0,1)$.
For $n=2$ and any $\bar\epsilon<2\epsilon<1$
$$\left\{\begin{array}{ll}
|D^jw(x)|\le C|x|^{2-\bar\epsilon-j},\quad |x|>2R_1,\quad j=0,1,2,3,4 \\
\\
\frac{|D^4w(y_1)-D^4w(y_2)|}{|y_1-y_2|^{\alpha}}\le C|y_1|^{-2-\bar \epsilon-\alpha},\quad |y_1|>2R_1,\,\, y_2\in B_{|y_1|/2}(y_1).
\end{array}
\right.
$$
\end{lem}

\noindent{\bf Proof of Lemma \ref{12feb27lem2}:}

Apply $\partial_k$ to $det(D^2 v)=f_v$ we have
\begin{equation}\label{12mar5e2}
a_{ij}\partial_{ij}(\partial_k v)=\partial_k f_v
\end{equation}
where $a_{ij}=cof_{ij}(D^2 v)$.
Lemma \ref{12feb27lem1} implies
$$|a_{ij}(x)-\delta_{ij}|\le \frac{C}{|x|^{\epsilon}},\quad |Da_{ij}(x)|\le \frac{C}{|x|^{1+\epsilon}},\quad |x|>R_1 $$
and for any $\alpha\in (0,1)$
$$\frac{|Da_{ij}(x_1)-Da_{ij}(x_2)|}{|x_1-x_2|^{\alpha}}\le C|x_1|^{-1-\epsilon-\alpha},\quad |x_1|>2R_1,
\quad x_2\in B_{|x_1|/2}(x_1). $$

Then apply $\partial_l$ to (\ref{12mar5e2}) and let $h_1=\partial_{kl}v$
$$
a_{ij}\partial_{ij}h_1=\partial_{kl}f_v-\partial_la_{ij}\partial_{ijk}v.
$$
We further write the equation above as
\begin{equation}\label{12feb29e1}
\Delta h_1=f_2:=\partial_{kl}f_v-\partial_la_{ij}\partial_{ijk}v-(a_{ij}-\delta_{ij})\partial_{ij}h_1.
\end{equation}

By (\ref{12feb23e1}) and Lemma \ref{12feb27lem1}, for any $\alpha\in (0,1)$
\begin{equation}\label{12feb22e5}
\left\{\begin{array}{ll}|f_2(x)|\le C|x|^{-2-2\epsilon}\quad |x|\ge 2R_1,\\
\\
\frac{|f_2(x_1)-f_2(x_2)|}{|x_1-x_2|^{\alpha}}\le \frac{C}{|x_1|^{2+2\epsilon+\alpha}},\,\,
x_2\in B_{|x_1|/2}(x_1),\,\, |x_1|\ge 2R_1.
\end{array}
\right.
\end{equation}
Note that by Lemma \ref{12feb27lem1} $h_1(x)\to \delta_{kl}$ as $x\to \infty$.
If $n\ge 3$, let
$$h_2(x)=-\int_{\mathbb R^n\setminus B_{R_1}}\frac{1}{n(n-2)\omega_n}|x-y|^{2-n}f_2(y)dy  $$
where $\omega_n$ is the volume of the unit ball in $\mathbb R^n$.
If $n=2$, let
$$h_2(x)=\frac{1}{2\pi}\int_{\mathbb R^2\setminus B_{R_1}}(\log |x-y|-\log |x|)f_2(y)dy. $$
In either case
$\Delta h_2=f_2$.
By elementary estimate it is easy to get
\begin{equation}\label{12feb27e5}
| D^jh_2(x)|\le \left\{\begin{array}{ll} C|x|^{-2\epsilon-j},\quad |x|>2R_1,\,\, j=0,1,\,\, n\ge 3, \\
C|x|^{-\bar \epsilon-j},\quad |x|>2R_1, \,\, j=0,1,\,\, n=2
\end{array}
\right.
\end{equation}
where $\bar \epsilon$ is any positive number less than $2\epsilon$.
Indeed, for each $x$, let
\begin{eqnarray*}
&&E_1:=\{ y\in \mathbb R^n\setminus B_{2R_1},\quad |y|\le |x|/2,\,\, \}, \\
&&E_2:=\{ y\in \mathbb R^n\setminus B_{2R_1},\quad |y-x|\le |x|/2,\, \},\\
&&E_3=(\mathbb R^n\setminus B_{2R_1})\setminus (E_1\cup E_2).
\end{eqnarray*}
Then it is easy to get (\ref{12feb27e5}).
For the estimate of $D^2h_2$ we claim that given $\alpha\in (0,1)$, if $n\ge 3$
\begin{equation}\label{12feb25e1}
\left\{\begin{array}{ll}
|D^jh_2(x)|\le C|x|^{-2\epsilon-j},\quad j=0,1,2,\quad |x|>2R_1, \\
\\
\frac{|D^2h_2(x_1)-D^2h_2(x_2)|}{|x_1-x_2|^{\alpha}}\le \frac{C}{|x_1|^{2+2\epsilon+\alpha}},\,\, x_2\in
B_{\frac{|x_1|}2}(x_1),\, |x_1|>2R_1.
\end{array}
\right.
\end{equation}
Replacing $2\epsilon$ by $\bar \epsilon$ we get the corresponding estimates of $D^2 h_2$ for $n=2$.
The way to obtain (\ref{12feb25e1}) is standard. Indeed,
 for each $x_0\in \mathbb R^n\setminus B_{2R_1}$, let $R=|x_0|$, we set
$$h_{2,R}(y)=h_2(x_0+\frac R4y),\quad f_{2,R}(y)=\frac{R^2}{16}f_2(x_0+\frac R4y), \quad |y|\le 2. $$
By (\ref{12feb22e5})
$\|f_{2,R}\|_{C^{\alpha}(B_1)}\le CR^{-2\epsilon}$.
Therefore Schauder estimate gives
$$\|h_{2,R}\|_{C^{2,\alpha}(B_1)}\le C(\|h_{2,R}\|_{L^{\infty}(B_2)}+\|f_{2,R}\|_{C^{\alpha}(B_2)})\le CR^{-2\epsilon},$$
which is equivalent to (\ref{12feb25e1}). The way to get the corresponding estimate for $n=2$ is the same.  Now we have
$$\Delta (h_1-h_2)=0,\quad \mathbb R^n\setminus B_{2R_1}. $$
Since we know $h_1-\delta_{kl}-h_2\to 0$ at infinity. For $n\ge 3$, by comparing with a multiple of $|x|^{2-n}$ we have
$$|h_1(x)-\delta_{kl}-h_2(x)|\le C|x|^{2-n},\quad |x|>2R_1.$$
By the estimate on $h_2$ we have
$$|h_1(x)-\delta_{kl}|\le C|x|^{-2\epsilon}, \quad |x|>2R_1. $$
Correspondingly
$$|D^jw(x)|\le C|x|^{2-j-2\epsilon},\quad |x|>2R_1,\quad j=0,1,2,\quad n\ge 3.$$
For $n=2$ we have
\begin{equation}\label{12mar7e1}
|h_1(x)-\delta_{kl}-h_2(x)|\le C|x|^{-1},\quad |x|>2R_1.
\end{equation}
Indeed, let $h_3(y)=h_1(\frac{y}{|y|^2})-\delta_{kl}-h_2(\frac{y}{|y|^2})$, then $\Delta h_3=0$ in $B_{1/2R_1}\setminus \{0\}$ and $\lim_{y\to 0}h_3(y)=0$. Therefore $|h_3(y)|\le C|y|$ near $0$. (\ref{12mar7e1}) follows. By fundamental theorem of calculus,
$$|D^jw(x)|\le C|x|^{2-j-\bar\epsilon},\quad |x|>2R_1,\quad j=0,1,2,\quad n=2.  $$
 Finally we
apply Lemma \ref{12feb27lem1} to obtain the estimates on the third and fourth derivatives.
Lemma \ref{12feb27lem2} is established. $\Box$

\medskip

{\bf Case one: $n\ge 3$. }

\medskip

Let $k_0$ be a positive integer such that $2^{k_0}\epsilon<1$ and $2^{k_0+1}\epsilon>1$ ( we choose $\epsilon$ smaller if necessary to make both inequalities hold). Let $\epsilon_1=2^{k_0}\epsilon$, clearly we have $1<2\epsilon_1<2$. Applying Lemma \ref{12feb27lem2} $k_0$ times we have
\begin{equation}\label{12mar7e3}
\left\{\begin{array}{ll}
|D^k w(x)|\le C|x|^{2-\epsilon_1-k},\quad k=0,..,4,\quad |x|>2R_1\\
\\
\frac{|D^4w(x_1)-D^4w(x_2)|}{|x_1-x_2|^{\alpha}}\le C|x_1|^{-2-\epsilon_1-\alpha},\,\, |x_1|>2R_1,\,\, x_2\in B_{|x_1|/2}(x_1).
\end{array}
\right.
\end{equation}
Let $h_1$ and $f_2$ be the same as in Lemma \ref{12feb27lem2}. Then we have
$$
\left\{\begin{array}{ll}|f_2(x)|\le C|x|^{-2-2\epsilon_1}\quad |x|\ge 2R_1,\\
\\
\frac{|f_2(x_1)-f_2(x_2)|}{|x_1-x_2|^{\alpha}}\le \frac{C}{|x_1|^{2+2\epsilon_1+\alpha}},\,\,
x_2\in B_{|x_1|/2}(x_1),\,\, |x_1|\ge 2R_1.
\end{array}
\right.
$$
Constructing $h_2$ as in Lemma \ref{12feb27lem2} ( the one for $n\ge 3$) we have
\begin{equation}\label{12mar11e2}
\left\{\begin{array}{ll}
|D^jh_2(x)|\le C|x|^{-2\epsilon_1-j},\quad j=0,1,2,\quad |x|>2R_1, \\
\\
\frac{|D^2h_2(x_1)-D^2h_2(x_2)|}{|x_1-x_2|^{\alpha}}\le \frac{C}{|x_1|^{2+2\epsilon_1+\alpha}},\,\, x_2\in
B_{\frac{|x_1|}2}(x_1),\,\, |x_1|>2R_1.
\end{array}
\right.
\end{equation}
As in the proof of Lemma \ref{12feb27lem2} by (\ref{12mar11e2}) we have
$$|h_1(x)-h_2(x)|\le C|x|^{2-n},\quad |x|>2R_1. $$
Since $2\epsilon_1>1$
$$|h_1(x)|\le |h_2(x)|+C|x|^{2-n}\le C|x|^{-1}. $$
By Theorem 4 of \cite{gilbarg}, $\partial_m w(x)\to c_m$ for some $c_m\in \mathbb R$ as $|x|\to \infty$.  Let $b\in \mathbb R^n$ be the limit of
$\nabla w$ and $w_1(x)=w(x)-b\cdot x$. The equation for $w_1$ can be written as ( for $e\in \mathbb S^{n-1}$)
\begin{equation}\label{12mar7e8}
a_{ij}\partial_{ij}(\partial_e w_1)=\partial_e f_v.
\end{equation}
By (\ref{12mar7e3}) the equation above can be written as
\begin{equation}\label{12mar7e4}
\Delta (\partial_e w_1)=f_3:=\partial_e f_v-(a_{ij}-\delta_{ij})\partial_{ije}w_1, \quad |x|>2R_1.
\end{equation}
and we have
$$\left\{\begin{array}{ll}
|f_3(x)|\le C(|x|^{-\beta-1}+|x|^{-1-2\epsilon_1})\le C|x|^{-1-2\epsilon_1},\,\, |x|>2R_1\\
\\
\frac{|f_3(x_1)-f_3(x_2)|}{|x_1-x_2|^{\alpha}}\le C|x_1|^{-1-2\epsilon_1-\alpha},\,\, |x_1|>2R_1,\,\, x_2\in B_{|x_1|/2}(x_1).
\end{array}
\right.
$$
Let $h_4$ solve $\Delta h_4=f_3$ and the construction of $h_4$ is similar to that of $h_2$. Then we have
$$\left\{\begin{array}{ll}
|D^jh_4(x)|\le C|x|^{1-2\epsilon_1-j},\,\, |x|>2R_1,\,\, j=0,1,2,\\
\\
\frac{|D^2h_4(x_1)-D^2h_4(x_2)|}{|x_1-x_2|^{\alpha}}\le C|x_1|^{-1-2\epsilon_1-\alpha},\,\, |x_1|>2R_1,\,\, x_2\in B_{|x_1|/2}(x_1).
\end{array}
\right.
$$
Since $\partial_e w_1-h_4\to 0$ at infinity, we have
\begin{equation}\label{12mar7e20}
|\partial_e w_1(x)-h_4(x)|\le C|x|^{2-n},\quad |x|>R_1.
\end{equation}
Therefore we have obtained $|\nabla w_1(x)|\le C|x|^{1-2\epsilon_1}$ on $|x|>R_1$. Using fundamental theorem of calculus
$$|w_1(x)|\le C|x|^{2-2\epsilon_1},\quad j=0,1, \quad |x|>R_1. $$
Lemma \ref{12feb27lem1} applied to $w_1$ gives
$$|D^jw_1(x)|\le C|x|^{2-j-2\epsilon_1}, \quad j=0..4. $$
Going back to (\ref{12mar7e4}), now the estimate for $f_3$ becomes
$$\left\{\begin{array}{ll}
|f_3(x)|\le C|x|^{-\beta-1}+C|x|^{-1-4\epsilon_1}, \quad |x|>2R_1, \\
\\
\frac{|f_3(x_1)-f_3(x_2)|}{|x_1-x_2|^{\alpha}}\le C(|x_1|^{-\beta-1-\alpha}+|x_1|^{-1-4\epsilon_1-\alpha}),\,\, |x_1|>2R_1, \,\,
x_2\in B_{|x_1|/2}(x_1).
\end{array}
\right.
$$
The new estimate of $h_4$ is
$$|h_4(x)|\le C(|x|^{1-\beta}+|x|^{1-4\epsilon_1}), \quad |x|>2R_1. $$
As before (\ref{12mar7e20}) holds.
Consequently
$$|\nabla w_1(x)|\le C(|x|^{2-n}+|x|^{1-4\epsilon_1})\le C|x|^{-1}, \quad |x|>2R_1. $$
By Theorem 4 of \cite{gilbarg}, $w_1\to c$ at infinity.  Let
$$w_2(x)=w(x)-b\cdot x-c. $$
Then we have $|w_2(x)|\le C$ for $|x|>2R_1$. Lemma \ref{12feb27lem1} applied to $w_2$ gives
\begin{equation}\label{12mar7e9}
 |D^k w_2(x)|\le C|x|^{-k},\quad k=0,1,2,3,\quad |x|>2R_1.
\end{equation}
The equation for $w_2$ can be written as
$$ det(I+D^2 w_2(x))=f_v. $$
Taking the difference between this equation and $det(I)=1$ we have
$$\tilde a_{ij}\partial_{ij}w_2=f_v-1,\quad |x|>2R_1 $$
where $\tilde a_{ij}$ satisfies
$$|D^j(\tilde a_{ij}(x)-\delta_{ij})|\le C|x|^{-2-j},\quad |x|>2R_1, \quad j=0,1. $$
Using (\ref{12mar7e9}) this equation can be written as
$$\Delta w_2=f_4:=f_v-1-(\tilde a_{ij}-\delta_{ij})\partial_{ij}w_2, \quad |x|>2R_1. $$
(\ref{12mar7e9}) further gives
$$\left\{\begin{array}{ll}
|f_4(x)|\le C(|x|^{-\beta}+|x|^{-4}),\quad |x|>2R_1, \\
\\
\frac{|f_4(x_1)-f_4(x_2)|}{|x_1-x_2|^{\alpha}}\le C(|x_1|^{-\beta-\alpha}+|x_1|^{-4-\alpha}),\quad |x_1|>2R_1, \,\,
x_2\in B_{|x_1|/2}(x_1).
\end{array}
\right.
$$
Let $h_5$ be defined similar to $h_2$. Then $h_5$ solves $\Delta h_5=f_4$ in $\mathbb R^n\setminus B_{2R_1}$ and satisfies
$$|h_5(x)|\le C(|x|^{2-\beta}+|x|^{-2}).$$
As before we have
$$|w_2(x)-h_5(x)|\le C|x|^{2-n}, \quad |x|>2R_1, $$
which gives
\begin{equation}\label{12mar11e3}
|w_2(x)|\le C(|x|^{2-n}+|x|^{2-\beta}+|x|^{-2}),\quad |x|>2R_1.
\end{equation}
If $|x|^{-2}>|x|^{2-n}+|x|^{2-\beta}$ we can apply the
same argument as above finite times to remove the $|x|^{-2}$ from (\ref{12mar11e3}). Eventually by Lemma \ref{12feb27lem1} we have (\ref{12mar7e10}). Proposition \ref{caf-li-prop} is established for $n\ge 3$.

\medskip

\noindent{\bf Case two: $n=2$}

As in the case for $n\ge 3$ we let $k_0$ be a positive integer such that $2^{k_0}\epsilon<1$ and $2^{k_0+1}\epsilon>1$ ( we choose $\epsilon$ smaller if necessary to make both inequalities hold). Let $\epsilon_1<2^{k_0}\epsilon$ and we let $1<2\epsilon_1<2$. Applying Lemma \ref{12feb27lem2} $k_0$ times then (\ref{12mar7e3}) holds.
Consider the equation for $w$. By taking the difference between the equation for $v$ and $det(I)=1$ we have
$$\tilde a_{ij}\partial_{ij}w=f_v-1. $$
We further write the equation above as
$$\Delta w=f_5:=f_v-1-(\tilde a_{ij}-\delta_{ij})\partial_{ij}w. $$
By (\ref{12mar7e3})
$$|f_5(x)|\le C|x|^{-2\epsilon_1},\quad |x|>R_1. $$
Let
$$h_6(x)=\frac{1}{2\pi}\int_{\mathbb R^2\setminus B_{R_1}}(\log |x-y|-\log |x|)f_5(y)dy. $$
Then elementary estimate gives
$$|h_6(x)|\le C|x|^{\epsilon_2}, \quad |x|>R_1 $$
for some $\epsilon_2\in (0,1)$. Since $w-h_6$ is harmonic on $\mathbb R^2\setminus B_{R_1}$ and $w-h_6=O(|x|^{2-\epsilon_1})$, there exist $b\in \mathbb R^2$ and $d_1,d_2\in \mathbb R$ such that
\begin{equation}\label{12mar8e1}
 w(x)-h_6(x)= b\cdot x+d_1\log |x|+d_2+O(1/|x|) \quad |x|>2R_1.
\end{equation}
Equation (\ref{12mar8e1}) is standard. For the convenience of the readers we include the proof. Let $z_l(r)$ be the projection of $w-h_6$ on $\sin l\theta$ for $l=1,2,..$. Then $z_l$ satisfies
$$z_l''(r)+\frac{1}{r}z_l'(r)-\frac{l^2}{r^2}z_l(r)=0,\quad r>2R_1. $$
Clearly $z_l(r)=c_{1l}r^l+c_{2l}r^{-l}$. Since $z_l(r)\le Cr^{2-\epsilon_1}$ we have $c_{1l}=0$ for all $l\ge 2$. Thus $z_l(r)=c_{2l}r^{-l}$. Let
$C$ be a constant such that $\max_{B_{2R_1}}|w-h_6|\le C$. Then $|z_l(2R_1)|\le C$, which gives $|c_{2l}|\le C(2R_1)^l$. The estimate for the projection of $w-h_6$ over $\cos l\theta$ for $l\ge 2$ is the same. The term $d_1\log |x|+d_2$ comes from the projection onto $1$. The projection onto $\cos\theta$ and $\sin\theta$ gives $b\cdot x$. (\ref{12mar8e1}) is established.

Let
$$w_1(x)=w(x)-b\cdot x. $$
Then $|w_1(x)|\le C|x|^{\epsilon_2}$. Apply Lemma \ref{12feb27lem1}
$$|D^kw_1(x)|\le C|x|^{\epsilon_2-k},\quad k=0,..,4, \quad |x|>2R_1. $$
The equation for $w_1$ can be written as
$$\Delta w_1=O(|x|^{-\beta})+O(|x|^{2\epsilon_2-4}). $$
Let
$$h_7(x)=\frac{1}{2\pi}\int_{\mathbb R^2\setminus B_{2R_1}}(\log |x-y|-\log |x|)\Delta w_1(y)dy. $$
Then
$$|h_7(x)|\le C(|x|^{2-\beta+\epsilon}+|x|^{2\epsilon_2-2+\epsilon}) $$
for $\epsilon>0$ arbitrarily small. Since $w_1-h_7$ is harmonic on $\mathbb R^2\setminus B_{2R_1}$ and
$w_1(x)-h_7(x)=O(|x|^{\epsilon_2})$, we have, for some $d,c\in \mathbb R$
$$w_1(x)-h_7(x)=d\log |x|+c+O(1/|x|). $$
Using the estimates on $h_7$ we have
\begin{equation}\label{12mar8e2}
w_1(x)=d\log |x|+c+O(|x|^{2\epsilon_2-2+\epsilon})+O(|x|^{2-\beta+\epsilon}).
\end{equation}
To obtain (\ref{12mar7e11}) we finally let
$$v_1(x)=v(x)-b\cdot x-c $$
and
$$H(x)=\frac 12|x|^2+d\log |x|. $$
Clearly $det(D^2 v_1(x))=f_v(x)$ and $det(D^2 H(x))=1-\frac{d^2}{|x|^4}$. Let $w_2(x)=v_1(x)-H(x)$. By (\ref{12mar8e2}) we already have
$$ |w_2(x)|\le C|x|^{-\epsilon_3},\quad |x|>2R_1 $$
for some $\epsilon_3>0$. Using Theorem \ref{caf-jian-wang} as well as Schauder estimate as in the proof of Lemma \ref{12feb27lem1} we obtain
\begin{equation}\label{12mar8e4}
|D^k w_2(x)|\le C|x|^{-\epsilon_3-k}\quad |x|>2R_1,\quad k=0,1,2.
\end{equation}
Thus the equation of $w_2$ can be written as
$$\Delta w_2(x)=O(|x|^{-4-2\epsilon_3})+O(|x|^{-\beta}). $$
Let
$$h_8(x)=\frac{1}{2\pi}\int_{\mathbb R^2\setminus B_{R_1}}(\log |x-y|-\log |x|)\Delta w_2(y)dy. $$
Then
\begin{equation}\label{12mar8e3}
|D^jh_8(x)|=O(|x|^{-2-j}+|x|^{\epsilon_5+2-\beta-j}),\quad j=0,1,\quad |x|>R_1
\end{equation}
for all $\epsilon_5>0$. Then we have $w_2(x)-h_8(x)=O(|x|^{-2})$ because of (\ref{12mar8e4}), (\ref{12mar8e3}) and the argument in the proof of (\ref{12mar8e1}). Consequently
$$w_2(x)=O(|x|^{-2}+|x|^{\epsilon_6+2-\beta}),\quad |x|>R_1 $$
for all $\epsilon_6>0$. The estimates on the derivatives of $w_2$ can be obtained by Lemma \ref{12feb27lem1}.
Proposition \ref{caf-li-prop} is established for $n=2$ as well. $\Box$

\section{Proof of Theorem \ref{thm5}}

Without loss of generality we assume that $B_2\subset D\subset B_{\bar r}$. First we prove a lemma that will be used in the proof for $n\ge 3$ and $n=2$.

\begin{lem} \label{thm5lem1}
There exists $c_1(n,\phi,D)$ such that for every $\xi\in \partial D$, there exists $w_{\xi}$ such that
$$\left\{\begin{array}{ll}
det(D^2w_{\xi}(x))\ge f(x)\quad \mathbb R^n\setminus D, \\
w_{\xi}(\xi)=\phi(\xi),\quad w_{\xi}(x)<\phi(x),\quad \forall x\in \partial D, \,\, x\neq \xi,
 \\
w_{\xi}(x)\le \frac 12|x|^2+c_1,\quad x\in (\mathbb R^n\setminus D)\cap B(0, 10 \mbox{diam}(D)).
\end{array}
\right.
$$
\end{lem}

\noindent{\bf Proof of Lemma \ref{thm5lem1}:} Let $f_1$ be a smooth radial function on $\mathbb R^n$ such that $f_1>f$ on $\mathbb R^n\setminus D$ and $f_1$ satisfies (FA). Let
$$z(x)=\int_0^{|x|}(\int_0^snt^{n-1}f_1(t)dt)^{\frac 1n}ds. $$
Then $det(D^2z(x))=f_1(x)$ on $\mathbb R^n$ and
$$|z(x)-\frac 12|x|^2|\le \left\{\begin{array}{ll} C,\quad n\ge 3, \\
C\log (2+|x|),\quad n=2
\end{array}
\right.
\quad x\in \mathbb R^n.  $$
Since $D$ is strictly convex, we can put $\xi$ as the origin using a translation and a rotation and then assume that $D$ stays in $\{x_n>0\}$. Assume that the boundary around $\xi$ is described by $x_n=\rho(x')$ where $x'=(x_1,..,x_{n-1})$. By the strict convexity we assume
$$\rho(x')=\frac 12 \sum_{1\le \alpha,\beta\le n-1}B_{\alpha\beta}x_{\alpha}x_{\beta}+o(|x'|^2) $$
where $(B_{\alpha\beta})\ge \delta I$ for some $\delta>0$.
By subtracting a linear function from $z$ we obtain $z_{\xi}$ that satisfies
$$\left\{\begin{array}{ll}
det(D^2z_{\xi})\ge f, \quad \mathbb R^n \setminus D, \\
z_{\xi}(0)=\phi(\xi), \quad \nabla z_{\xi}(0)=\nabla \phi(\xi)
\end{array}
\right.
$$
and
$$|z_{\xi}(x)-\frac 12|x|^2|\le C|x|,\quad x\in \mathbb R^n\setminus D. $$
Next we further adjust $z_{\xi}$ by defining
$$w_{\xi}(x)=z_{\xi}(x)-A_{\xi}x_n $$
for $A_{\xi}$ large to be determined.
When evaluated on $\partial D$ near $0$,
$$w_{\xi}(x',\rho(x'))-\phi(x',\rho(x'))\le C|x'|^2-A_{\xi}\rho(x'). $$
Therefore for $|x'|\le \delta_1$ for some $\delta_1$ small we have $w_{\xi}(x',\rho(x'))<\phi(x',\rho(x'))$. For $|x'|>\delta_1$, the convexity of $\partial D$ yields
$$x_n\ge \delta_1^3, \quad \forall x\in \partial D\setminus \{(x',\rho(x')):\quad |x'|<\delta_1. \}. $$
Then by choosing $A_{\xi}$ possibly larger (but still under control) we have $w_{\xi}(x)<\phi(x)$ for all $x\in \partial D$. Clearly $A_{\xi}$ has a uniform bound for all $\xi\in \partial \Omega$.
Lemma \ref{thm5lem1} is established. $\Box$

\medskip

Let
$$\underline{w}(x)=\max\left\{w_{\xi}(x)~\big|~\xi\in\partial{D}\right\}.$$
It is clear by Lemma \ref{thm5lem1} that $\underline{w}$ is a
locally Lipschitz function in $B_{2\bar r}\setminus{D}$, and
$\underline{w}=\varphi$ on $\partial{D}$. Since $w_{\xi}$ is a
smooth convex solution of \eqref{extdir}, $\underline{w}$ is a
viscosity subsolution of \eqref{extdir} in
$B_{2\bar r}\setminus\overline{D}$.  Let $c_1$ be the constant determined in Lemma \ref{thm5lem1}. Then we have
$$\underline{w}(x)\le \frac 12|x|^2+c_1,\quad B_{2\bar r}\setminus \bar D. $$

\medskip

We finish the proof of Theorem \ref{thm5} in two cases.

\noindent{\bf Case one: $n\ge 3$.}  Clearly we only need to prove the existence of solutions for $A=I$ and $b=0$, as the general case can be reduced to this case by a linear transformation. Let
$\bar f$ and $\underline{f}$ be smooth, radial functions such that
$\underline{f}<f<\bar f$ in $\mathbb R^n\setminus D$ and suppose $\underline{f}$ and $\bar f$
satisfy (FA). For $d>0$ and $\beta_{1},\beta_{2}\in\mathbb{R}$, set
$$\underline{u}_{d}(x)=\beta_{1}+\int_{\bar r}^r\bigg (\int_1^snt^{n-1}\bar f(t)dt+d\bigg )^{\frac 1n}ds, \quad r=|x|>2,$$
and
$$\overline{u}_{d}(x)=\beta_{2}+\int_{2}^r\bigg (\int_1^snt^{n-1}\underline{f}(t)dt+d\bigg )^{\frac 1n}ds, \quad r=|x|>2.$$
Clearly
$$\det(D^2\underline{u}_{d})=\bar f\ge f,\quad \mathbb R^n\setminus \bar D,$$
and
$$\det(D^2\overline{u}_{d})=\underline{f}\leq f,\quad \mathbb R^n\setminus \bar D.$$
On the other hand,
\begin{equation}\label{h1-1}
\underline{u}_{d}(x)\leq\beta_{1},\quad\quad\mbox{in}~B_{\bar{r}}\setminus\overline{D},~\forall~d>0.
\end{equation}
and
\begin{equation}\label{h1-12}
\overline{u}_{d}(x)\geq\beta_{2},\quad\quad\mbox{in}~B_{\bar{r}}\setminus\overline{D},~\forall~d>0.
\end{equation}
Let
$$ \beta_{1}:=\min\left\{\underline{w}(x)~\big| ~x\in\overline{B_{\bar{r}}}\setminus{D}\right\}-1<\min_{\partial{D}}\varphi,$$
$$ \beta_{2}:=\max_{\partial{D}}\varphi+1.$$
This shows that
$\underline{u}_{d}$ and $\overline{u}_{d}$ are continuous convex
subsolution and supersolution of \eqref{12apr5e1}, respectively. By the definition of $\underline{u}_d$ by choosing $d$ large enough, say $d\ge d_0$, we can make
$$\underline{u}_d> \underline{w}(x)+1,\quad |x|=\bar r+1. $$
By (\ref{h1-1}) and the above, the function
$$u_{1,d}(x)=\left\{\begin{array}{ll} \underline{u}_d,\quad |x|\ge \bar r+1, \\
\underline{w}(x),\quad x\in B_{\bar r}\setminus D,\\
\max\{ \underline{w}(x), \underline{u}_d\},\quad x\in B_{\bar r+1}\setminus B_{\bar r}
\end{array}
\right.
$$
is a viscosity subsolution of (\ref{12apr5e1}) if $d\ge d_0$.

Next we consider the asymptotic behavior of $\underline{u}_d$ and $\bar u_d$ when $d$ is fixed.
Using (FA) it is easy to obtain
$$\underline{u}_{d}(x)=\frac 12|x|^2+\mu_{1}(d)+O(|x|^{2-\min\{\beta,n\}}), $$
and
$$\overline{u}_{d}(x)=\frac 12|x|^2+\mu_{2}(d)+O(|x|^{2-\min\{\beta,n\}}), $$
where
$$\mu_{1}(d)=\beta_{1}-\frac{\bar{r}^{2}}{2}+\int_{\bar{r}}^{\infty}\left(\bigg (\int_1^snt^{n-1}\bar f(t)dt+d\bigg )^{\frac 1n}-s\right)ds,$$
and
$$\mu_{2}(d)=\beta_{2}-2+\int_{2}^{\infty}\left(\bigg (\int_1^snt^{n-1}\underline{f}(t)dt+d\bigg )^{\frac 1n}-s\right)ds.$$
It is easy to see that $\mu_{1}(d)$ and $\mu_{2}(d)$ are strictly
increasing functions of $d$ and
\begin{equation}\label{mud}
\lim_{d\to\infty}\mu_{1}(d)=\infty,\quad\mbox{and}\quad\lim_{d\to\infty}\mu_{2}(d)=\infty.
\end{equation}

Let $c_*=\mu_1(d_0)$, recall that for $d>d_0$, $u_{1,d}$ is a viscosity subsolution.
For every $c>c_{*}$, there exists a unique
$d(c)$ such that
\begin{equation}\label{h4-1}
\mu_{1}(d(c))=c.
\end{equation}
So $\underline{u}_{d(c)}$ satisfies
\begin{equation}\label{h5-1}
\underline{u}_{d(c)}(x)=\frac{1}{2}|x|^{2}+c+O\left(|x|^{2-\min\{\beta,n\}}\right),
\quad\quad\mbox{as}\quad\,x\rightarrow\infty.
\end{equation}

Also there exists $d_{2}(c)$ such that
$\mu_{2}(d_{2}(c))=c$ and
\begin{equation}\label{h5-12}
\overline{u}_{d_{2}(c)}(x)=\frac{1}{2}|x|^{2}+c+O\left(|x|^{2-\min\{\beta,n\}}\right),
\quad\quad\mbox{as}\quad\,x\rightarrow\infty.
\end{equation}
By \eqref{h5-1} and \eqref{h5-12}
$$\lim_{|x|\rightarrow\infty}\left(\underline{u}_{d(c)}(x)-\overline{u}_{d_2(c)}(x)\right)=0.$$
On the other hand, by the definition of $\beta_1$ we have $\bar u_{d_2(c)}>u_{1,d(c)}$ on $\partial D$.
Thus, in view of the comparison principle for smooth convex
solutions of Monge-Amp\`{e}re, (see \cite{cns}), we have
\begin{equation}\label{h10-1}
u_{1,d(c)}\leq \bar u_{d_2(c)},\quad\mbox{on}~\mathbb{R}^{n}\setminus{D}.
\end{equation}

For any $c>c_{*}$, let $\mathcal{S}_{c}$ denote the set of
$v\in{C}^{0}(\mathbb{R}^{n}\setminus{D})$ which are viscosity
subsolutions of \eqref{12apr5e1} in
$\mathbb{R}^{n}\setminus\overline{D}$ satisfying
\begin{equation}\label{h13-1}
v=\varphi,\quad\mbox{on}~\partial{D},
\end{equation}
and
\begin{equation}\label{h13-2}
u_{1,d(c)}\leq{v}\leq \bar u_{d_2(c)},\quad\mbox{in}~\mathbb{R}^{n}\setminus{D}.
\end{equation}
We know that $u_{1,d(c)}\in\mathcal{S}_{c}$. Let
$$u(x):=\sup\left\{v(x)~|~v\in\mathcal{S}_{c}\right\}, \quad x \in \mathbb{R}^{n}\setminus{D}.$$
Then $u$ is convex and of class
${C}^{0}(\mathbb{R}^{n}\setminus{D})$. By \eqref{h5-1}, and the
definitions of $u_{1,d(c)}$ and $\bar u_{d_2(c)}$
\begin{equation}\label{h13-3}
u(x)\geq u_{1,d(c)}(x)=\frac{1}{2}|x|^{2}+c+O\left(|x|^{2-\min\{\beta,n\}}\right),
\quad\quad\mbox{as}\quad\,x\rightarrow\infty
\end{equation}
and
$$ u(x)\leq \overline{u}_{d_2(c)}(x)=\frac{1}{2}|x|^{2}+c+O\left(|x|^{2-\min\{\beta,n\}}\right). $$
The estimate \eqref{11may5e1} for $k=0$ follows.

Next, we prove that $u$ satisfies the boundary condition. It is
obvious from the definition of $u_{1,d(c)}$ that
$$ \liminf_{x\rightarrow\xi}u(x)\geq \lim_{x\rightarrow\xi}u_{1,d(c)}(x)
=\varphi(\xi),\quad\forall~\xi\in\partial{D}.$$ So we only need to
prove that
$$\limsup_{x\rightarrow\xi}u(x)\leq \varphi(\xi),\quad\forall~\xi\in\partial{D}.$$
Let $\omega_{c}^{+}\in C^2(\overline{B_{\bar{r}}\setminus D})$ be
defined by
$$
\begin{cases}
\Delta\omega_{c}^{+}=0,&\mbox{in}~B_{\bar{r}+1}\setminus{\overline{D}},\\
\omega_{c}^{+}=\varphi,&\mbox{on}~\partial{D},\\
\omega_{c}^{+}=\max\limits_{\partial{B_{\bar{r}+1}}}\overline{u}_{d_2(c)},&\mbox{on}~\partial{B_{\bar{r}+1}}.
\end{cases}
$$
It is easy to see that a viscosity subsolution $v$ of \eqref{12apr5e1}
satisfies $\Delta{v}\geq0$ in viscosity sense. Therefore, for every
$v\in\mathcal{S}_{c}$, by $v\leq\omega_{c}^{+}$ on
$\partial(B_{\bar{r}}\setminus{D})$, we have
$$v\leq\omega_{c}^{+}\quad\mbox{in}~B_{\bar{r}}\setminus{\overline{D}}.$$
It follows that
$$u\leq\omega_{c}^{+}\quad\mbox{in}~B_{\bar{r}}\setminus{\overline{D}},$$
and then
$$\limsup_{x\rightarrow\xi}u(x)\leq \lim_{x\rightarrow\xi}\omega_{c}^{+}(x)
=\varphi(\xi),\quad\forall~\xi\in\partial{D}.$$

Finally, we prove $u$ is a solution of \eqref{extdir}. For
$\bar{x}\in\mathbb{R}^{n}\setminus\overline{D}$, fix some $\epsilon>0$
such that
$B_{\epsilon}(\bar{x})\subset\mathbb{R}^{n}\setminus\overline{D}$.
By the definition of $u$, $u\leq\bar{u}$. We claim that there is a convex viscosity solution to
$\widetilde{u}\in{C}^{0}(\overline{B_{\epsilon}(\bar{x})})$ to
$$
\begin{cases}
\det(D^{2}\widetilde{u})=f,&x\in{B}_{\epsilon}(\bar x),\\
\widetilde{u}=u,&x\in\partial{B}_{\epsilon}(\bar x).
\end{cases}
$$
Indeed, let $\phi_k$ be a sequence of smooth functions on $\partial B_{\epsilon}(\bar x)$ satisfying
$$u\le \phi_k\le u+\frac 1k. $$
Let $f_k$ be a sequence of smooth positive functions tending to $f$ and $f_k\le f$. Let $\psi_k$ be the convex solution to
$$\left\{\begin{array}{ll}
det(D^2\psi_k)=f_k\quad B_{\epsilon}(\bar x),\\
\psi_k=\phi_i\quad \mbox{ on }\quad \partial B_{\epsilon}(\bar x).
\end{array}
\right.
$$
Clearly $\psi_k\ge u$. On the other hand, let $h_k$ be the harmonic function on $B_{\epsilon}(\bar x)$ with $h_k=\phi_k$ on $\partial B_{\epsilon}(\bar x)$. Then we have $u_k \le h_k$. Therefore $|\psi_k|$ is uniformly bounded over any compact subset of $B_{\epsilon}(\bar x)$. $|\nabla \psi_k|$ is also uniformly bounded over all compact subsets of $B_{\epsilon}(\bar x)$ by the convexity. Thus $\psi_k$ converges along a subsequence to $\tilde u$ in $B_{\epsilon}(\bar x)$. By the closeness between $h_k$ to $u$ on $\partial B_{\epsilon}(\bar x)$, $\tilde u$ can be extended as a continuous function to $\bar B_{\epsilon}(\bar x)$.
By the maximum principle, $u\leq\widetilde{u}\leq\bar{u}_{d_2(c)}$ on
${B}_{\epsilon}$. Define
$$w(y)=
\begin{cases}
\widetilde{u}(y),&\mbox{if}~y\in{B}_{\epsilon},\\
u(y),&\mbox{if}~y\in\mathbb{R}^{2}\setminus(D\cup{B}_{\epsilon}(\bar{x})).
\end{cases}
$$
Clearly, $w\in\mathcal{S}_{c}$. So, by the definition of $u$,
$u\geq{w}$ on ${B}_{\epsilon}(\bar{x})$. It follows that
$u\equiv\widetilde{u}$ on ${B}_{\epsilon}(\bar{x})$. Therefore $u$
is a viscosity solution of \eqref{extdir}. We have proved (\ref{11may5e1}) for $k=0$. The estimates of derivatives follow from Proposition \ref{caf-li-prop}. Theorem \ref{thm5} is established for $n\ge 3$.

\medskip

\noindent{\bf Case two: $n=2$}.

As in case one we let $\bar f$ be a radial function such that $\bar f(|x|)\ge f(x)$ in $\mathbb R^2\setminus D$, and $\bar f$ also satisfies (FA). Let
$$\underline{u}_d(x)=\beta_1+\int_{\bar r}^r\bigg (\int_1^s 2t\bar f(t)dt+d\bigg )^{\frac 12}ds $$
for $d\ge 0$ and $r>1$. Here we choose $\beta_1=\min_{\partial D}\phi-1$. Clearly
 $$\underline{u}_d(x)<\underline{w}(x)\quad B_{\bar r}\setminus D,\quad \forall d\ge 0. $$
Then we choose $d^*$ large so that for all $d\ge d^*$, $\underline{u}_d(x)>\underline{w}(x)$ on $\partial B_{\bar r+1}$. Let
$$u_{1,d}(x)=\left\{\begin{array}{ll}
\underline{w}(x),\quad B_{\bar r}\setminus \bar D\\
\max\{\underline{w}(x),\underline{u}_d\},\quad B_{\bar r+1}\setminus \bar B_{\bar r},\\
\underline{u}_d,\quad \mathbb R^2\setminus B_{\bar r+1}.
\end{array}
\right.
$$
Then $u_{1,d}$ is a convex viscosity subsolution of (\ref{12apr5e1}).  Let
$$A_d=d-1+\int_1^{\infty}2t(\bar f(t)-1)dt. $$
Then elementary computation gives
$$\underline{u}_d(x)=\frac 12|x|^2+A_d\log |x|+O(1). $$
Next we let $\underline{f}$ be a radial function such that $\underline{f}(|x|)\le f(x)$ for $x\in \mathbb R^2\setminus \bar D$. Suppose $\underline{f}$ also satisfies (FA) and is positive and smooth on $\mathbb R^2$. Let
$$\bar u_d(x)=\beta_2+\int_2^r\bigg (\int_1^s2t\underline{f}(t)dt+d\bigg )^{\frac 12}ds. $$
Let
$$L_d=d-1+\int_1^{\infty}2t(\underline{f}(t)-1)dt. $$
Then the asymptotic behavior of $\bar u_d$ at infinity is
$$\bar u_d(x)=\frac 12|x|^2+L_d\log |x|+O(1). $$
Thus for all $d>d^*$, we can choose $d_1$ such that $L_{d_1}=A_d$.
Then we choose $\beta_1$ such that $\bar u_{d_1}>\phi$ on $\partial D$ and $\bar u_{d_1}>\underline{u}_d$ at infinity.
As in case one, by taking the supremum of subsolutions we obtain a solution $u$ that is equal to $\phi$ on $\partial D$ and
$$u(x)=\frac 12|x|^2+A_d\log |x|+O(1). $$
By Proposition \ref{caf-li-prop}
$$u(x)=\frac 12|x|^2+A_d\log |x|+c+o(1). $$
The following lemma says the constant term is uniquely determined by other parameters.

\begin{lem}\label{uniquenesslem}
Let $u_1$, $u_2$ be two locally convex smooth functions on $\mathbb R^2\setminus \bar D$ where $D$ satisfies the same assumption as in Theorem \ref{thm5}. Suppose
$u_1$ and $u_2$ both satisfy
$$\left\{\begin{array}{ll}det(D^2u)=f\mbox{ in }\quad \mathbb R^2\setminus \bar D,\\
u=\phi,\quad \mbox{ on }\quad \partial D
\end{array}
\right.
$$ with $f$ satisfying (FA) and for the same constant $d$
\begin{equation}\label{12feb8e1}
u_i(x)-\frac 12|x|^2-d\log |x|=O(1), \quad x\in \mathbb R^2\setminus \bar D,\quad i=1,2.
\end{equation}
Then $u_1\equiv u_2$.
\end{lem}

\noindent{\bf Proof of Lemma \ref{uniquenesslem}:} By Proposition \ref{caf-li-prop} we see that when (\ref{12feb8e1}) holds, we have
$$D^2u_i(x)=I+O(|x|^{-2+\epsilon}), \quad i=1,2 $$
for $\epsilon>0$ small and  $|x|$ large. For the proof of this lemma we only need
\begin{equation}\label{11dec23e1}
D^2u_i(x)=I+O(|x|^{-\frac 32}),\quad |x|>1,\quad i=1,2.
\end{equation}
By Proposition \ref{caf-li-prop},
$$u_i(x)=\frac 12|x|^2+d\log |x|+c_i+O(1/|x|^{\sigma}), \quad i=1,2 $$
for $\sigma\in (0, \min\{\beta-2,2\})$.
Without loss of generality we assume $c_1>c_2$. If $c_1=c_2$ we know $u_1\equiv u_2$ by maximum principle. Since $u_1=u_2$ on $\partial D$, we have,
$u_1>u_2$ in $\mathbb R^2\setminus \bar D$. Let $w=u_1-u_2$, then $w$ satisfies
$$a_{ij}\partial_{ij}w=0, \quad \mathbb R^2\setminus \bar D $$
where
$$a_{ij}(x)=\int_0^1cof_{ij}(tD^2u_1+(1-t)D^2u_2)dt. $$
By the assumption of Lemma \ref{uniquenesslem} and (\ref{11dec23e1}), $a_{ij}$ is uniformly elliptic and
\begin{equation}\label{11dec23e2}
a_{ij}(x)=\delta_{ij}+O(|x|^{-\frac 32}),\quad x\in \mathbb R^2\setminus \bar D.
\end{equation}
Let $a_0<\frac 12a_1$ be positive constants to be determined. We set
$h_{\epsilon}=\epsilon \log (|x|-a_0)$ over $a_1<|x|<\infty$. Direct computation shows, by (\ref{11dec23e2}) that
\begin{eqnarray*}
&&a_{ij}\partial_{ij}h_{\epsilon}=\Delta h_{\epsilon}+(a_{ij}-\delta_{ij})\partial_{ij}h_{\epsilon}\\
&\le & -\frac{\epsilon a_0}{(|x|-a_0)^2|x|}+C\epsilon |x|^{-7/2}, \\
&\le & -\frac{4\epsilon a_0}{|x|^3}+C\epsilon |x|^{-7/2}, \quad |x|>a_1>a_0.
\end{eqnarray*}
By choosing $a_0$ sufficiently large and $a_1>2a_0$ we have
$$a_{ij}\partial_{ij}h_{\epsilon}<0, \quad a_1<|x|<\infty. $$
Let $R>a_1$ and $M_R=\max_{|x|=R}w$. Let $v=w-M_R$, then clearly for all $\epsilon>0$, $h_{\epsilon}$ is greater than $v$ on $\partial B_R$ and at infinity. Thus for any compact subset $K\subset\subset \mathbb R^2\setminus \bar B_R$, $v<h_{\epsilon}$. Let $\epsilon\to 0$ we have
$$w(x)\le M_R, \quad \forall |x|\ge R. $$
Taking any $R_1>R$, we have $\max_{|x|=R_1}w\le M_R$. Strong maximum principle implies that either $\max_{|x|=R_1}w<M_R$ for all $R_1>R$ or $w$ is a constant. $w$ is not a constant, therefore we have $\max_{|x|=R_1}w<M_R$ for all $R_1>R$. However, this means over the region $B_{R_1}\setminus \bar D$, the maximum of $w$ is attained at an interior point, a contradiction to the elliptic equation that $w$ satisfies. Thus Lemma \ref{uniquenesslem} is established. $\Box$

\medskip

Lemma \ref{uniquenesslem} uniquely determines the constant in the expansion, then by Proposition \ref{caf-li-prop} we obtain (\ref{12jan18e1}). Thus Theorem \ref{thm5} for the case $n=2$ is established. $\Box$

\section{Proof of Theorem \ref{thm3}}

We only need to consider the existence part as the uniqueness part follows immediately from maximum principles. For the existence part we only need to consider the case that $A=I$, $b=0$ and $c=0$, because the general case can be reduced to this case by a linear transformation. Consider $u_R$ that solves
\begin{equation}\label{12jan19e1}
\left\{\begin{array}{ll}
det(D^2u_R)=f, \quad B_R, \\
\\
u_R=\frac{R^2}2, \quad \partial B_R.
\end{array}
\right.
\end{equation}
We shall bound $u_R$ above and below by two radial functions.
Let $h$ be a smooth radial function, then at the point $(|x|, 0,...,0)$
$$D^2h(x)=\mbox{diag}(h''(r), h'(r)/r,...,h'(r)/r), \quad r=|x|. $$
Thus $det(D^2h)(x)=h''(r)(h'(r)/r)^{n-1}. $

We first construct a subsolution $h_-(r)$: Let
$\bar f$ be a radial function such that $\bar f>f$ and $\bar f$ satisfies (FA).
$$h_-(r)=\int_0^r(\int_0^s nt^{n-1}\bar f(t)dt)^{\frac 1n}ds. $$
Clearly $det(D^2h_-)=\bar f$ in $\mathbb R^n$ and since
$\bar f(t)=1+O(t^{-\beta})$ it is
easy to verify that
$$h_-(r)=\frac 12|x|^2+O(1). $$
Next we construct a super solution. Let
$\underline f$ be a radial function less than $f(x)$ and $\underline f$ also satisfy (FA),
$$h_+(r)=\int_0^r(\int_0^s n t^{n-1}\underline f(t)dt)^{\frac 1n}ds. $$
Similarly we have $det(D^2h_+)=\underline{f}$ in $\mathbb R^n$ and
$h_+(r)=\frac 12r^2+O(1)$ for $r$ large.
Let $\beta_-$ be a constant such that $h_-(|x|)+\beta_-\le \frac 12|x|^2$, $\beta_+$ be a constant such that $h_+(|x|)+\beta_+\ge \frac 12|x|^2$. Then by maximum principle
\begin{equation}\label{11dec29e1}
h_-(r)+\beta_-\le u_R(x)\le h_+(r)+\beta_+, \quad |x|\le R.
\end{equation}
Let $R\to \infty$ and the sequence $u_R$ converges to a global solution $u$ that satisfies $det(D^2u)=f$ in $\mathbb R^n$ and $u-\frac 12|x|^2=O(1)$. For this convergence, we use the fact that for any $K\subset\subset \mathbb R^n$, $|u_R(x)-\frac 12|x|^2|\le C(K)$ and
by Caffarelli's $C^{1,\alpha}$ estimate \cite{caf-cpam-ca}, $\|\nabla u_R\|_{L^{\infty}(K)}\le C(K)$. Thus $u_R$ converges to a convex viscosity solution $u$ to $det(D^2u)=f$ in $\mathbb R^n$ with the property that
$$|u(x)-\frac 12|x|^2|\le C, \quad \mathbb R^n. $$
By Proposition \ref{caf-li-prop}, there exists a $c^*\in \mathbb R$ such that
$$\lim_{|x|\to \infty} |x|^{\min\{\beta,n\}-2+k}\bigg (D^k(u(x)-\frac 12|x|^2-c^*)\bigg )<\infty $$
for $k=0,1,2,3,4$.
After a translation the solution with the desired asymptotic behavior can be found. Theorem \ref{thm3} is established. $\Box$

\section{The proof of Theorem \ref{thm1}}

Without loss of generality we assume $u(0)=0=\min_{\mathbb R^n}u$. The goal is to show that there exists a linear transformation $T$ such that
$v=u\cdot T$ satisfies (\ref{12feb27e2}). Then we employ Proposition \ref{caf-li-prop} to finish the proof. The proof of $v$ satisfying (\ref{12feb27e2}) is by the argument of Caffarelli-Li.

Suppose $c_0^{-1}\le \inf_{\mathbb R^n}f\le \sup_{\mathbb R^n}f<c_0$, only under this assumption it is proved in \cite{caf-li1} that for $M$ large and
$$\Omega_M:=\{x\in \mathbb R^n;\quad u(x)<M \quad \}$$
there exists $a_M\in \mathcal{A}$ such that
\begin{equation}\label{11may5e2}
B_{R/C}\subset a_M(\Omega_M)\subset B_{CR},
\end{equation}
where $R=\sqrt{M}$ and $C>1$ is a constant independent of $M$.
Let
$$O:=\{y;\quad a_M^{-1}(Ry)\in \Omega_M \}. $$
Then $B_{1/C}\subset O\subset B_C$. Set
$$\xi(y):=\frac{1}{R^2}u(a_M^{-1}(Ry)),$$
then we have
\begin{equation}\label{eqforxi}
\left\{\begin{array}{ll}
det(D^2\xi)=f(a_M^{-1}(Ry)),\quad \mbox{in}\quad O, \\
\xi=1,\quad \mbox{on}\quad \partial O.
\end{array}
\right.
\end{equation}
Let $\bar \xi$ solve
$$\left\{\begin{array}{ll}
det(D^2\bar \xi)=1,\quad \mbox{in}\quad O, \\
\bar \xi=1,\quad \mbox{on}\quad \partial O.
\end{array}
\right.
$$
By Pogorelov's estimate
$$ \frac 1C I\le D^2\bar \xi \le CI, \quad |D^3\bar \xi(x)|\le C, \,\, x\in O,\,\, dist(x,\partial O)\ge \delta. $$
We claim that there exists $C>0$ independent of $M$ such that
\begin{equation}\label{11may5e3}
|\xi(x)-\bar \xi(x)|\le C/R, \quad x\in O.
\end{equation}
Indeed, by the Alexandrov estimate (\cite{caf-cab})
$$-\min_{\bar O}(\xi-\bar \xi)\le C\bigg (\int_{S^+}det(D^2(\xi-\bar \xi))\bigg )^{1/n} $$
where
$$S^+:=\{x\in O;\quad D^2(\xi-\bar \xi)>0 \quad \}. $$
On $S^+$
$$\frac{D^2\xi}2=\frac{D^2(\xi-\bar \xi)}2+\frac{D^2\bar \xi}2, $$
so the concavity of $det^{\frac 1n}$ on positive definite symmetric matrices implies
$$det(D^2(\xi-\bar \xi))^{\frac 1n}\le f(a_M^{-1}(Ry))^{\frac 1n}-1. $$
Thus
$$-\min_{\bar O}(\xi-\bar \xi)\le C\bigg (\int_{S^+}|f(a_M^{-1}(Ry))^{\frac 1n}-1|^ndy)^{\frac 1n}. $$
Let $z=a_M^{-1}(Ry)$, i.e. $a_Mz=Ry$ then $dz=R^ndy$
$$
\bigg(\int_{S^+}|f(a_M^{-1}(Ry))^{\frac 1n}-1)^n|dy\bigg )^{\frac 1n}
\le \frac 1R\bigg (\int_{B_{CR}}|f(z)^{\frac 1n}-1|^ndz\bigg )^{\frac 1n}. $$
By the assumption (FA) the integral is finite, thus
we have proved that
$$-\min_{\bar O}(\xi-\bar \xi)\le C/R,\quad x\in O $$
Similarly we also have
$-\min_{\bar O}(\bar \xi-\xi)\le C/R.$
(\ref{11may5e3}) is proved.

Next we set
$$E_M:=\{x;\quad (x-\bar x)'D^2\xi(\bar x)(x-\bar x)\le 1 \quad \}$$
where $\bar x$ is the minimum of $\bar \xi$. By Theorem 1 of \cite{caf-ann1} $\bar x$ is the unique minimum point of $\bar \xi$. Then by the same argument as in \cite{caf-li1} we have the following: There exist $\bar k$ and $C$ depending only on $n$ and $f$ such that for $\epsilon=\frac 1{10}$,
$M=2^{(1+\epsilon)k}$, $2^{k-1}\le M'\le 2^k$,
$$(\frac{2M'}{R^2}-C2^{-\frac{3\epsilon k}{2}})^{\frac 12}E_M\subset \frac 1Ra_M(\Omega_{M'})\subset
(\frac{2M'}{R^2}+C2^{-\frac{3\epsilon k}{2}})^{\frac 12}E_M, \quad \forall k\ge \bar k, $$
which is
$$\sqrt{2M'}(1-\frac{C}{2^{\epsilon k/2}})E_M\subset a_M(\Omega_{M'})\subset \sqrt{2M'}(1+\frac{C}{2^{\epsilon k/2}})E_M. $$
Let $Q$ be a positive definite matrix such that $Q^2=D^2\bar \xi(\bar x)$, $O$ be an orthogonal matrix such that
$T_k:=O Q_k a_M$ is upper triangular. Then clearly $det(T_k)=1$ and by Proposition 3.4 of \cite{caf-li1} we have
$$\|T_k-T\|\le C2^{-\frac{\epsilon k}2}. $$
Let $v=u\cdot T$, then clearly
$$det(D^2v(x))=f(Tx). $$
For $v$ and some $\bar k$ large we have
\begin{eqnarray*}
\sqrt{2M'}(1-\frac{C}{2^{\epsilon k/2}})B\subset \{x;\quad v(x)<M'\}\\
\subset \sqrt{2M'}(1+\frac{C}{2^{\epsilon k/2}})B \quad \forall M'\ge 2^{\bar k}.
\end{eqnarray*}
Consequently
\begin{equation}\label{11may9e1}
|v(x)-\frac 12 |x|^2|\le C|x|^{2-\epsilon}.
\end{equation}

Clearly $f(T\cdot )$ also satisfies (FA). Proposition \ref{caf-li-prop} gives the asymptotic behavior of $u$ and the estimates on its derivatives. the constant $d$ in the estimate in two dimensional spaces is determined similarly as in \cite{caf-li1}. Theorem \ref{thm1} is established .
$\Box$

\medskip

\begin{rem}
Corollary \ref{cor1} follows from Theorem \ref{thm1} just like in \cite{caf-li1} so we omit the proof.
\end{rem}

\section{Appendix: Interior estimate of Caffarelli and Jian-Wang}

The following theorem is a combination of the interior estimate of Caffarelli \cite{caf-ann2} and an improvement by Jian-Wang \cite{jian-wang}.

\medskip

\begin{thm} \label{caf-jian-wang} (Caffarelli, Jian-Wang) Let $u\in C^0(\Omega)$ be a convex viscosity solution of
\begin{eqnarray*}
&& det(D^2u)=f, \quad \Omega , \\
&& u=0\quad \mbox{on }\quad \partial \Omega,
\end{eqnarray*}
where $\Omega$ is a convex bounded domain satisfying $B_1\subset \Omega\subset B_n$. Assume that $f$ is Dini continuous on $\Omega$ and
$$ \frac 1{c_0}\le f \le c_0, \quad \Omega. $$
Then $u\in C^2(B_{1/2})$ and $\forall x,y\in B_{1/2}$
\begin{equation}\label{12march28e1}
|D^2u(x)-D^2u(y)|\le C\big ( d+\int_0^d\frac{\omega_f(r)}r+d\int_d^1\frac{\omega_f(r)}{r^2}\big )
\end{equation}
where $d=|x-y|$, $C>0$ depends only on $n$ and $c_0$, $\omega_f$ is the oscillation function of $f$ defined by
$$\omega_f(r):=\sup\{|f(x)-f(y)|:\quad |x-y|\le r \}. $$
 It follows that
(i) If $f$ is Dini continuous, then $u\in C^2(B_{1/2})$, and the modulus of convexity of $D^2u$ can be estimated by (\ref{12march28e1}).
(ii) If $f\in C^{\alpha}(\Omega)$ and $\alpha\in (0,1)$, then
$$\|D^2u\|_{C^{\alpha}(B_{1/2})}\le C\big(1+\frac{\|f\|_{C^{\alpha}(\Omega)}}{\alpha(1-\alpha)}\big ). $$
(iii) If $f\in C^{0,1}(\Omega)$, then
$$|D^2u(x)-D^2u(y)|\le Cd\big (1+\|f\|_{C^{0,1}(\Omega)}|\log d|\big ). $$
\end{thm}

 Here we recall that $f$ is Dini continuous if the oscillation function $\omega_f$ satisfies $\int_0^1\omega_f(r)/rdr<\infty$. 
 
 \begin{rem} Note that in 
 Caffarelli's interior estimate $u=0$ is assumed on $\partial \Omega$. Since $\Omega$ is very close to a ball, by \cite{caf-ann1,caf-cpam-ca} $u$ is strictly convex in $\Omega$. But there is no explicit formula that describes how the higher order derivatives of $u$ depend on l$f$. In Jian-Wang's theorem, this dependence is given as in (\ref{12march28e1}) but instead of assume $u=0$ on $\partial \Omega$, they assumed $u$ is strictly convex and their constant depends on the strict convexity. We feel the way that Theorem \ref{caf-jian-wang} is stated is more convenience for application.
 We only used the $(ii)$ and $(iii)$ of Theorem \ref{caf-jian-wang} in this article. 
 \end{rem}

\end{document}